\input amstex
\input amsppt.sty
\nologo
\parskip 5pt minus 2pt
\baselineskip=14pt plus 3pt
\hoffset.4truein\voffset-.125truein
\hsize5.75truein\vsize8.5truein
\loadeusm
\loadbold


\redefine\B{{\text{\rm B}}}
\redefine\C{{\Bbb C}}

\redefine\N{{\Bbb N}}
\redefine\O{{\text{\rm O}}}

\redefine\Q{{\Bbb Q}}
\redefine\R{{\Bbb R}}
\redefine\Z{{\Bbb Z}}

\define\br#1{{\left<{#1}\right>}}
\define\dd#1#2{{\frac{\partial{#1}}{\partial{#2}}}}
\define\df{{\,\overset\text{\rm df}\to =}\,}
\define\edge#1{{\,\overline{\!#1}}}
\define\endpf{{\hfill $\blacksquare$}}

\redefine\mod{{\text{\rm mod\,}}}

\define\w#1{{\widehat{#1}}}

\define\Card{{\text{\rm Card}}}
\define\Pf{{\flushpar{\it Proof.\;\;}}}

\flushpar
{\bf
Fixed Points of the $\bold p$-Adic $\bold q$-Bracket
\newline
Eric Brussel
\newline
Emory University}

\

\flushpar
{\bf Abstract.}
The $q$-bracket $[X]_q:\O_{\C_p}\to\O_{\C_p}$, which is the $q$-analog of the
identity function, is also a norm-preserving isometry, for each $q\in \B(1,p^{-1/(p-1)})$.
In this paper we investigate its fixed points.

\

\flushpar
{\bf I. Introduction.}

We start with an example from complex analysis.
Let $D$ be the unit disk in the complex plane $\C$.
An {\it isometry} of $D$ is a continuous,
distance-preserving map from $D$ to $D$.
All analytic isometries of $D$ are rotations, and preserve the complex norm.
They are parameterized in a natural way by $\R/\Z$,
with $t\in\R/\Z$ corresponding to the rotation $\rho_t:z\mapsto ze^{2\pi i t}$.
The quotient topology on $\R/\Z$ makes the isometries
into a {\it continuous family}, since for all $z\in D$,
$\lim_{t\to t_0}\rho_t(z)=\rho_{t_0}(z)$.
The fixed point set of this family is uninteresting,
since a nontrivial rotation fixes only the origin.
A more interesting set of fixed points is provided by
the larger family of {\it analytic automorphisms} of $D$.
By Schwarz's Lemma, this family is continuously parameterized by $\R/\Z\times D$,
with $(t,z_0)$ corresponding to the mobius transformation
$z\mapsto e^{2\pi it}\frac{z_0-z}{1-\bar z_0 z}$.
A direct computation shows
that an analytic automorphism has either one interior fixed point,
or one boundary fixed point, or two boundary fixed points.

This paper grew from an interest in the fixed points 
on the {\it $p$-adic} unit disk $\Z_p$.
Let $p$ be a prime, and let $\Z_p$ denote the additive group of $p$-adic integers.
There is a continuous family of norm-preserving isometries
$$
[X]_q:\Z_p\;\longrightarrow\;\Z_p
$$
parameterized by the elements $q$ of
the topological group $\B(1,p^{-1/(p-1)})$, which is $1+p\Z_p$ if $p$ is odd, and $1+4\Z_2$ if $p=2$.
The function $[X]_q$ is called the {\it $q$-bracket},
and it is an interpolation to $\Z_p$ of the arithmetic function
on $\N\cup\{0\}$ given by
$$
[n]_q=1+q+q^2+\cdots+q^{n-1}.
$$
The $q$-bracket is also known as the {\it $q$-analog} (or {\it $q$-extension})
of the identity function, and its values are {\it $q$-numbers}.
It is the canonical 1-cocycle $[X]_q\in\text{Z}^1(\Z_p,\Z_p)$ sending $1$ to $1$, 
where $\Z_p$ is viewed as a $\Z_p$-module via the action $1*1=q$.  Since $q$ is in $\B(1,p^{-1/(p-1)})$,
we have
$[X]_q\in\text{Z}^1(\Z_p,\Z_p)$ and $[X]_q(\mod p^n)\in\text{Z}^1(\Z/p^n,\Z/p^n)$.

Our results show that if $p\neq 3$, or if $q\equiv 1(\mod p^2)$, then $[X]_q$ has only the
``trivial'' fixed points $0$ and $1$ in $\Z_p$. 
However, if $p=3$, $q\equiv 1(\mod 3)$, and $q\not\equiv 1(\mod 9)$,
then $[X]_q$ has a unique nontrivial fixed point in $\Z_3$ for all $q$.
For example, if $p=3$ and $q=4$, then $-1/2$ is the nontrivial fixed point of $[X]_q$ when $q=4$: 
$[-1/2]_4=-1/2$.
The admissible $q$ in $\B(1,3^{-1/2})$ form two (disjoint) balls, $\B(4,3^{-1})$ and $\B(7,3^{-1})$,
the set of nontrivial fixed points $x$ for these $q$
respectively form the two balls $\B(1,1)$ and $\B(0,1)$,
and the map $Q(X)$ taking $x$ to $q$ is a bijective analytic contraction by the factor $3$.

It turns out to be easier to analyze our problem in the following more general context.

Let $\C_p$ denote the $p$-adic complex numbers.
Write $|-|$ for the metric on $\C_p$,
and let $\O_p=\O_{\C_p}$ denote the unit disk in $\C_p$.
Write $v$ for the corresponding additive valuation,
so that $|x|=p^{-v(x)}$.
Let $\B(a,r)$ denote the set $\{b\in\C_p:|b-a|<r\}$. 
If $x\in \O_p$ and $q\in \B(1,p^{-1/(p-1)})$ then 
$x\cdot\log q\in x\cdot\B(0,p^{-1/(p-1)})\subset \B(0,p^{-1/(p-1)})$, so
$q^x=\exp(x\log q)$
is well defined.  We define the {\it $q$-bracket} $[X]_q$ on $\O_p$ by
$$
[x]_q\df\;\cases
\frac{q^x-1}{q-1} &\text{ if $q\neq 1$ }\\
x&\text{ if $q=1$ }
\endcases
$$
As above,
the $q$-bracket is an interpolation to $\O_p$ of the $q$-number function on $\N$, defined by
$[n]_q=1+q+q^2+\cdots+q^{n-1}$.
The $q$-bracket is also the canonical 1-cocycle $[X]_q\in\text{Z}^1(\O_p,\O_p)$ sending $1$ to $1$, 
where $\O_p$ is viewed as an $\O_p$-module via the action $1*1=q$.  

In this setting, we will show that 
$x\in\O_{\C_p}$ is a nontrivial fixed point of $[X]_q$ for some $q$
if and only if $|A_{p-2}(x)|>p^{-1/(p-1)}$, where $A_{p-2}(X)\df(X-2)\cdots(X-(p-1))$,
and then $|A_{p-2}(x)|=|p(q-1)^{2-p}|$.
The set of nontrivial pairs $(x,q)$ such that $[x]_q=x$ form a manifold whose
standard projections each have degree $p-2$.
If $[x]_q=x$, then we have a surjective analytic map $Q(X):\B(x,|A_{p-2}(x)|)\to\B(q,|q-1|)$ such that
$[x']_{Q(x')}=x'$ for all $x'\in\B(x,|A_{p-2}(x)|)$, 
and this map is a (bijective) contraction if and only if the residue $\bar x$ has multiplicity one in the fiber over $q$,
if and only if $\bar x$ is not a root of the polynomial $A_{p-2}^{(1)}(X)$, i.e., 
if and only if $|A_{p-2}^{(1)}(x)|=1$.

For reference on basic concepts see the beautiful book \cite{G} by Gouv\^ea.
The study of $q$-functions for a general variable $q$ tending to 1 is old,
and the study of $q$-numbers and $q$-identities goes back at least to Jackson in \cite{J}.
In \cite{F} Fray proved $p$-adic $q$-analogs of theorems of Legendre, Kummer, and
Lucas on $q$-binomial coefficients.
The structure of the space of continuous functions
$C(K,\Q_p)$, where $K$ is a local field, was studied by Dieudonn\'e in \cite{D},
and Mahler constructed an explicit basis for this space in \cite{M}.
In \cite{C} Conrad proved that the set of $q$-binomial coefficients (which includes the $q$-identity function),
form a basis for $C(\Z_p,\Z_p)$.
Isometries on $\Z_p$ or on locally compact connected one-dimensional abelian groups were studied in
\cite{A}, \cite{B}, and \cite{Su}.

We would like to thank the referee of a previous version of this paper for suggestions that dramatically
simplified the proofs and improved the results.

\

\flushpar
{\bf II. Results.}

\proclaim{Notation}
\rm
The letters $q,u$, and $x$ will always
denote elements of $\B(1,p^{-1/(p-1)}),\O_p^*$, and $\O_p$, respectively,
and the first two 
will frequently be related by $q=1+p^{m_0}u$, where $m_0=v(q-1)>1/(p-1)$.
The capitals $Q,U$, and $X$ will denote coordinate functions defined on these sets,
related to each other in a similar manner.
\endproclaim

\proclaim{Proposition 1}
Fix $q\in \B(1,p^{-1/(p-1)})$.
Then $[X]_q:\O_p\to \O_p$ is a norm-preserving isometry.
\endproclaim

\Pf
This is clear for $q=1$, so assume $q\neq 1$.
$[X]_q$ is the composition of analytic isomorphisms
$$
\O_p\,@>\;\log q\;>>\; \B(0,|q-1|)\,@>\;\exp(-)-1\;>>\;\B(0,|q-1|)\,@>\frac 1{q-1}\;>>\;\O_p
$$
a dilation by $|q-1|$, an isometry, and a contraction by $|q-1|$.
Tracing through the maps shows $[X]_q$ preserves the norm. 

\endpf

Since $[X]_q$ is an isometry of the $p$-adic unit disk onto itself, 
the notion of fixed point makes sense.
Set 
$$
f(X,Q)=[X]_Q-X
$$
For fixed $q\in \B(1,p^{-1/(p-1)})$ the set of fixed points of $[X]_q$ is the set of solutions 
$\{x:f(x,q)=0\}$.
We see that $f(X,Q)$ is analytic on $\O_p\times \B(1,p^{-1/(p-1)})$, since 
$$
f(X,Q)=\sum_{n=1}\binom X{n+1}(Q-1)^n
$$
is in $\C_p[[X,Q-1]]$ and converges on $\O_p\times \B(0,p^{-1/(p-1)})$.

It is obvious that $f(x,1)=0$ for all $x\in \O_p$, and $f(0,q)=f(1,q)=0$ for all $q\in \B(1,p^{-1/(p-1)})$.
We call these solutions {\it trivial}, and define
the set of {\it nontrivial fixed points}
$$
M\;\df\;V\left(\frac{[X]_Q-X}{(Q-1)X(X-1)}\right)
\;=\;\left\{(x,q):[x]_q=x\right\}
$$

\proclaim
{Proposition 2}
The set $M$ is a submanifold of $\O_p\times \B(1,p^{-1/(p-1)})$.
If $(x,q)\in M$
then there is an analytic function $Q(X)$ in a neighborhood $N$ of
$x$ such that $q=Q(x)$, and $(x',Q(x'))\in M$ for all $x'\in N$.
\endproclaim

\Pf
Define
$$
g(X,Q)\;\df\;\frac{f(X,Q)}{(Q-1)X(X-1)}
$$
We show $dg$ does not vanish on $M$ by showing $\dd gQ(x,q)\neq 0$ or all $(x,q)\in M$.
Then $M$ is a submanifold of $\O_p\times \B(1,p^{-1/(p-1)})$ by \cite{Se, Chapter III, Section 11},
and there is an analytic function $Q(x)$ such that
$(x',Q(x'))\in M$ in a neighborhood of $(x,q)$ by the $p$-adic implicit function theorem.
Directly from the definition of $f$,
$$
\dd fQ=\dd{[X]_Q}Q=\frac{XQ^{X-1}-[X]_Q}{Q-1}
$$
If $(x,q)\in M$, then $(x,q)\in V(f)$, so $[x]_q=x$.  
Therefore 
$$
\dd fQ(x,q)=x[x-1]_q
$$
Since $g=f/(Q-1)X(X-1)$,
$$
\dd gQ=\frac 1{(Q-1)X(X-1)}\dd fQ-\frac{g}{Q-1}
$$
By the power series expression for $g$ we have $g(X,1)=1/2$,
so in particular if $(x,q)\in M$ then $q\neq 1$, hence
$(g/(Q-1))(x,q)=0$, hence
$$
\dd gQ(x,q)=\frac 1{(q-1)x(x-1)}\dd fQ(x,q)=\frac{[x-1]_q}{(q-1)(x-1)}
$$
Using $q^X= e^{X\log q}$, we compute
$$
[X]_q=\sum_{n=1}^\infty\frac{(\log q)^n}{(q-1)n!}X^n
$$
Therefore $\frac{[X-1]_q}{X-1}=\sum_{n=0}\frac{(\log q)^{n+1}}{(q-1)(n+1)!}(X-1)^n$.
This has (nonzero) value $\log q/(q-1)$ at $X=1$, and if $X\neq 1$ then it is nonzero since
$[X-1]_q$ preserves the norm in $\O_p$ by Proposition 1.
We conclude $\dd gQ(x,q)\neq 0$ for all $(x,q)\in M$.

\endpf

\proclaim{Definitions}\rm
Let $\phi_1:M\to \O_p$ and $\phi_2:M\to \B(1,p^{-1/(p-1)})$
be the projections, and
let $M_q=\phi_2^{-1}(q)$ denote the fiber of $\phi_2$ over $q$.
We identify $M_q$ with $\phi_1(M_q)$, which is the set of nontrivial fixed points
of $[X]_q$.
\endproclaim

\flushpar
{\bf Series A.}
For any $x\in \O_p$ and $q\in \B(1,p^{-1/(p-1)})$, 
$$
\align
[X]_q-X
&=\sum_{n=0}^\infty c_n(x,q)(X-x)^n
\\
&=
([x]_q-x)+\left(\frac{q^x\log q}{q-1}-1\right)(X-x)+\sum_{n=2}
\frac{q^x(\log q)^n}{(q-1)n!}(X-x)^n
\endalign
$$
This series converges on $\O_p$.
Note that since $c_2$ is nonzero, any fixed point can have a maximum multiplicity of two.

\proclaim
{Proposition 3}
If $p=2$ then $M=\varnothing$.
If $p\neq 2$, then $\phi_2$ has degree $p-2$, and 
$$
\phi_2(M)
=\B(1,p^{-1/(p-1)})-\B(1,p^{-1/(p-2)})
$$
\endproclaim

\Pf
Fix $q\in \B(1,p^{-1/(p-1)})$, and let $m_0=v(q-1)$.  For $x\in \O_p$, let $c_n=c_n(x,q)$ be the coefficient from Series A. 
By the $p$-adic Weierstrass preparation theorem 
\cite{G, Theorem 6.2.6}, the number of zeros of $[X]_q-X$ (in $\O_p$) is $N=\sup\{n:v(c_n)=\inf_m v(c_m)\}$,
counting multiplicities.
Since $\{0,1\}$ are both zeros, we know $N\geq 2$, and since $M$ excludes these solutions,
$M_q$ has cardinality $N-2$.
We compute $N$:
If $n\geq 2$, 
$$
v(c_n)
=(n-1)m_0-\frac{n-s_p(n)}{p-1}
$$
where $s_p(n)$ is the sum of the coefficients of the $p$-adic expansion of $n$.
It is easy to see $v(c_n)>v(c_p)$ whenever $n>p$, and
$$
v(c_n(x,q))=\cases 
(n-1)m_0 &\text{ if $2\leq n\leq p-1$}\\
(p-1)m_0-1 &\text{ if $n=p$}
\endcases
\tag $*$
$$
Thus if $p=2$ or $m_0>1/(p-2)$ then $N=2$, hence $M_q=\varnothing$.
If $p\neq 2$ and $m_0\leq 1/(p-2)$ then $v(c_p)\leq v(c_2)$, hence
$N=p$, hence $M_q$ has cardinality $p-2$, counting multiplicities.
We conclude $\phi_2$ has degree $p-2$ for $1/(p-1)<m_0\leq 1/(p-2)$.

\endpf

Set $A_0(X)=1$, and for $n>0$, set
$$
A_n(X)\;\df\;(X-2)(X-3)\cdots(X-(n+1))
$$
Let $A_n^{(i)}(X)$ denote the $i$-th (formal) derivative.

\flushpar
{\bf Series B.}
Set $U=p^{-m_0}(Q-1)$ for $m_0>1/(p-1)$.  For any $u\in\O_p$ we compute
$$
g(x,Q)=\frac{[x]_Q-x}{(Q-1)x(x-1)}=\sum_{n=0}^\infty d_n(x,u)(U-u)^n
$$
where $d_n(x,u)=\sum_{k=n}^\infty\binom kn \frac{A_k(x)}{(k+2)!}p^{km_0} u^{k-n}$.
Note $d_n(x,0)=\frac{A_n(x)}{(n+2)!}p^{nm_0}$.

\proclaim{Proposition 4}
Suppose $p\neq 2$.
Then $\phi_1$ has degree $p-2$, and
$$\phi_1(M)=\{x\in\O_p:|A_{p-2}(x)|>p^{-1/(p-1)}\}$$
If $(x,q)\in M$, then $|A_{p-2}(x)|=p^{(p-2)m_0-1}$, where $m_0=v(q-1)$.
\endproclaim

\Pf
We use Series B with $u=0$.  Set $d_n=d_n(x,0)$.
Then 
$$
v(d_n)=v(A_n(x))+nm_0+\frac{s_p(n+2)-(n+2)}{p-1}
$$
and from this we read off
$$
v(d_n(x,0))=\cases
v(A_n(x))+nm_0&\text{ if $0\leq n\leq p-3$}\\
v(A_{p-2}(x))+(p-2)m_0-1 &\text{ if $n=p-2$}
\endcases
\tag $**$
$$ 
If this series has a solution $U=u\in\O_p^*$, then 
$p\neq 2$ and $1/(p-2)\geq m_0>1/(p-1)$ by Proposition 3.
If $n>p-2$ and $m_0>1/(p-1)$, then using the fact that $A_{p-2}(X)$ divides
$A_n(X)$, we easily compute $v(d_n)-v(d_{p-2})>0$:
$$
\align
v(d_n)-v(d_{p-2})
&=v(A_n(x))-v(A_{p-2}(x))+(n-(p-2))m_0+1+\frac{s_p(n+2)-(n+2)}{p-1}
\\
&>v(A_n(x))-v(A_{p-2}(x))+\frac{s_p(n+2)-1}{p-1}
\geq 0
\endalign
$$
Thus the Weierstrass polynomial has nonzero degree if and only if $v(d_{p-2})\leq v(d_0)$, i.e.,
$v(A_{p-2}(x))\leq 1-(p-2)m_0$, in which case the degree is $p-2$.
For a given $x$ this holds for some $m_0$ in the range $1/(p-2)\geq m_0>1/(p-1)$ if and only if
$v(A_{p-2}(x))<1/(p-1)$.
Now given $x\in\O_p$ such that $|A_{p-2}(x)|>p^{-1/(p-1)}$,
set $m_0=(1-A_{p-2}(x))/(p-2)$.  Then $v(d_0)=v(d_{p-2})$, so that all $p-2$
solutions are units $u$ such that $(x,q)\in M$, where $q=1+p^{m_0}u$.

\endpf

\proclaim{Proposition 5}
For each $q\in\B(1,p^{-1/(p-1)})-\B(1,p^{-1/(p-2)})$, let $m_0=v(q-1)$, and
let $\edge M_q$ denote the set of residues of the elements $\phi_1(\phi_2^{-1}(q))$.
\roster
\item"{a)}"
If $m_0<1/(p-2)$, $\edge M_q=\{2,\dots,p-1\}$, and $\Card(M_q)=p-2$.
\item"{b)}"
If $m_0=1/(p-2)$, $\edge M_q\cap\{2,\dots,p-1\}=\varnothing$, and $\Card(M_q)\geq p-3$. 
\endroster
In particular, $m_0=1/(p-2)$ if and only if $|A_{p-2}(x)|=1$.
\endproclaim

\Pf
The Weierstrass polynomial for $[X]_q-X$ in Series A 
has degree $p$ by Proposition 3.  
Suppose $x\in M_q$.
Then $v(A_{p-2}(x))=0$ if and only if 
$\bar x\not\in\{2,\dots,p-1\}$, if and only if $m_0=1/(p-2)$ by Proposition 4.  
Therefore $\edge M_q\subset\{2,\dots,p-1\}$ if $m_0<1/(p-2)$,
and $\edge M_q\cap\{2,\dots,p-1\}=\varnothing$ if $m_0=1/(p-2)$.

If $m_0<1/(p-2)$ then $v(c_n(x,q))>v(c_p(x,q))$ for $n=2,\dots,p-1$ by $(*)$.
Since not every fixed point of $[X]_q$ has the same residue we must have $v(c_1(x,q))=v(c_p(x,q))$,
by \cite{G, Corollary 6.4.11},
hence there is at most one $x\in M_q$ with any given residue, and $\Card(M_q)=p-2$.

If $m_0=1/(p-2)$ then $v(c_2(x,q))=v(c_p(x,q))=m_0$ by $(*)$, 
and the Newton polygon for Series A shows
there are at most two zeros with residue $\bar x$, using $(*)$ and \cite{G, Corollary 6.4.11}.
Suppose $[X]_q$ has fixed points $x,x'$, and $x''$, such that $x\neq x''$ and $\bar x=\bar x''\neq\bar x'$.
We compute $c_1(x',q)-c_1(x,q)=([x']_q-[x]_q)\log q$.  Since $x'-x$ has nonzero residue it is a unit,
hence $[x']_q-[x]_q$ is a unit by Proposition 1, so $v(c_1(x',q)-c_1(x,q))=m_0$.
Since $x\neq x''$ and $\bar x=\bar x''$, we have
$v(c_1(x,q))>m_0$ by the Newton polygon, and it follows that $v(c_1(x',q))=m_0$, 
so that $x'$ is the only fixed point with residue $\bar x'$.
Thus there are at most two points in $M_q$ with the same residue, 
hence $\Card(M_q)\geq p-3$.
The last statement is immediate.

\endpf

\flushpar
{\bf Remark.}
If $0$ or $1$ is in $M_q$ then $|M_q|=p-2$, since these are also trivial fixed points.
By Proposition 4 and Proposition 5, we compute
$$
\phi_1(M)=
\bigcup_{\underset{\bar a\not\in\{\bar 2,\dots,\bar p-\bar 1\}}\to{a\in\O_p}}\B(a,1)\;\cup\;\bigcup_{a\in\{2,\dots,p-1\}}
\B(a,1)-\overline\B(a,p^{-1/(p-1)})
$$
where the left union corresponds to $v(q-1)=1/(p-2)$, the right union to $v(q-1)<1/(p-2)$.
Note no rational integer not congruent to $0$ or $1$ $(\mod p)$ may be a fixed point of any $[X]_q$.

\proclaim{Theorem 6}
Suppose $x\in\O_p$ 
is a nontrivial fixed point
of $[X]_q$, for some $q\in\B(1,p^{-1/(p-1)})$.
Then $\Card(\phi_1^{-1}(x))=p-2$, and
each $(x,q)\in\phi_1^{-1}(x)$ determines a distinct residue $\bar u$,
for $q=1+p^{-v(q-1)}u$.
If $(x,q)\in M$ and $x'\in \B(x,|A_{p-2}(x)|)$ then there exists a unique $q'\in \B(q,|q-1|)$ such that
$(x',q')\in M$, and the resulting map
$$
Q(X):\B(x,|A_{p-2}(x)|)\to \B(q,|q-1|)
$$
defined by $Q(x')=q'$ is an analytic surjection satisfying $|q'-q|=|([x']_q-x')/x'(x'-1)|$.
The map is a (bijective) contraction (by $p^{1-(p-1)v(q-1)}$) if and only if $|A_{p-2}^{(1)}(x)|=1$,
if and only if the multiplicity of $\bar x$ in $\edge M_q$ is one.
This occurs for all but finitely many residue classes for $x$.
\endproclaim

\Pf
Let $m_0=v(q-1)$.
Since $x$ is a nontrivial fixed point, $\phi_1$ has degree $p-2$ by Proposition 4.
We show the cardinality of $\phi_1^{-1}(x)$ is $p-2$ by showing
that the various $u$ appearing in $q=1+p^{m_0}u\in\phi_1^{-1}(x)$ have distinct
residues.
Suppose $(x,q)\in\phi_1^{-1}(x)$, $q=1+p^{m_0}u$,
and $d_n(x,u)$ is the coefficient of Series B for $g(x,1+p^{m_0}U)$, expanded around $u$.
Then $d_0(x,u)=0$,
and using $(**)$, the identity $d_n(x,u)=\sum_{i=0}\binom{n+i}n d_{n+i}(x,0)u^i$,
and the fact that $v(d_{p-2}(x,0))=0$ (by Proposition 4),
we compute $v(d_n(x,u)) =0$ for $1\leq n\leq p-2$.
Thus the Newton polygon contains the points 
$(n,v(d_n(x,u)))=(0,\infty),(1,0),\dots,(p-2,0)$,
so that no other solution $U=u'$ has the same residue as $u$,
by \cite{G, Corollary 6.4.11}.
We conclude $\Card(\phi_1^{-1}(x))=p-2$, and the $p-2$ roots $u$ have distinct residues.

Next, for $x'\in\B(x,|A_{p-2}(x)|)$,
we compute $|g(x',q)|$ in terms of $|x'-x|$.
We will show that $|g(x',q)|=p^{1-(p-2)m_0}|x'-x|$ except when ($m_0=1/(p-2)$
and) $|A_{p-2}^{(1)}(x)|=1$, in which case $|g(x',q)|<|x'-x|$, 
and that $|A_{p-2}^{(1)}(x)|=1$ for finitely many residue classes $\bar x$.

Since $g(x,q)=0$, by Series B we have
$$
g(x',q)=(x'-x)\sum_{k=1}a_k=(x'-x)\sum_{k=1}
\frac{D_k(x')}{(k+2)!}p^{km_0}u^k
$$
where $D_k(X)=(A_k(X)-A_k(x))/(X-x)$.
Compute 
$$
\align
v(a_k)&=v(D_k(x'))+km_0-\frac{k+2-s_p(k+2)}{p-1}\\
&= v(D_k(x'))+k\delta+\frac{s_p(k+2)-2}{p-1}
\endalign
$$
where $\delta=m_0-1/(p-1)>0$.
If $k\neq p^r-2$ then $s_p(k+2)-2\geq 0$, so $v(a_k)>0$ in this case.
If $k=p^r-2$, $v(a_k)=v(D_{p^r-2}(x'))+(p^r-2)\delta-1/(p-1)$.

If $m_0<1/(p-2)$ then $|A_{p-2}(x)|<1$ by Proposition 5,
and $|A_{p-2}^{(1)}(x)|=1$
since $A_{p-2}(X)$ is separable $(\mod p)$.
But $|A_{p^r-2}^{(1)}(x)|<1$ for $r>1$, since $x$ is a multiple root of $A_{p^r-2}(X)(\mod p)$
for $r>1$.
Since
$$
D_k(x')=
A_k^{(1)}(x)+\frac 1{2!}A_k^{(2)}(x)(x'-x)+\cdots+\frac 1{k!}A_k^{(k)}(x)(x'-x)^{k-1}
$$
we conclude $v(D_{p^r-2}(x'))>v(D_{p-2}(x'))=0>v(A_{p-2}(x))-v(x'-x)$ 
for $r>1$ and $x'\in\B(x,|A_{p-2}(x)|)$. 
Since $(p^r-2)\delta-(p-2)\delta=(p^r-p)\delta>0$ for $r>1$,
$v(a_{p^r-2})>v(a_{p-2})$ for $m_0<1/(p-2)$, and so $v(a_{p-2})=(p-2)m_0-1$ 
is the unique minimum value of all of the coefficients.  
We conclude $|g(x',q)|=p^{1-(p-2)m_0}|x'-x|=|A_{p-2}(x)|^{-1}|x'-x|$,
and since $x'\in\B(x,|A_{p-2}(x)|)$, this shows
$1>|g(x',q)|>|x'-x|$ when $m_0<1/(p-2)$.

If $m_0=1/(p-2)$, then $|A_{p-2}(x)|=1$ by Proposition 5.
We compute 
$$
A_{p-2}^{(1)}(X)=(p-2)X^{p-3}+(p-3)X^{p-4}+\cdots+2X+1(\mod p)
$$
Thus for all but finitely many residue classes $\bar x$, we have $A_{p-2}^{(1)}(x)\neq 0(\mod p)$,
hence $D_{p-2}(x')=A_{p-2}^{(1)}(x)\neq 0(\mod p)$ whenever $x'\in\B(x,1)$, hence $v(a_{p-2})=0$.
We also compute $v(a_k)\geq m_0$ for $k\neq p^r-2$ for some $r$, and 
since $v(a_{p^r-2})\geq (p^r-p)/(p-1)(p-2)>0$ for $r>1$, this shows
$|g(x',q)|=|x'-x|$ when ($m_0=1/(p-2)$ and) $|A_{p-2}^{(1)}(x)|=1$.

If $m_0=1/(p-2)$ and $|A_{p-2}^{(1)}(x)|<1$, then $|D_{p-2}(x')|<1$,
so $v(a_{p-2})>(p-2)m_0-1=0$.  
Checking by hand, we find $A_{p-2}^{(1)}(0)=A_{p-2}^{(1)}(1)=1(\mod p)$,
so $\bar x\neq 0,1$, and $A_{p-2}(1/2)=0(\mod p)$, so $\bar x\neq 1/2$.
Since $A_{p^r-2}(X)=A_{p-2}(X)^r(X(X-1))^{r-1}(\mod p)$,
the product rule
shows that when $r>1$, $A_{p^r-2}^{(1)}(x)\neq 0(\mod p)$, since $\bar x\neq 0,1,1/2$.
Therefore $|D_{p^r-2}(x')|=1$, hence 
$v(a_{p^r-2})=(p^r-p)/(p-2)(p-1)>0$ for $r>1$.
We conclude at any rate that $|g(x',q)|<|x'-x|$ when ($m_0=1/(p-2)$ and) $|A_{p-2}^{(1)}(x)|<1$.

We have shown that $|A_{p-2}^{(1)}(x)|=1$ implies $|g(x',q)|=p^{1-(p-2)m_0}|x'-x|$,
hence that $g(x',q)\neq 0$ if $\bar x'=\bar x$, i.e., the
multiplicity of $\bar x$ in $\edge M_q$ is one.
Conversely, suppose $x\in M_q$ and $\bar x$ has multiplicity one in $\edge M_q$.
We have already seen that $|A_{p-2}^{(1)}(x)|=1$
if $\bar x\in\{0,1\}$, and otherwise
$v(c_1(x,q))=v(c_p(x,q))(=(p-1)m_0-1)$ are the minimum values
in Series A for $f(X,q)$, and
$v(f(x',q))=v(c_1(x,q)(x'-x))=v(x'-x)+(p-1)m_0-1$.
Therefore $|g(x',q)|=p^{1-(p-2)m_0}|x'-x|\geq |x'-x|$, for all $x'\in\B(x,1)$,
and it follows from the above that $|A_{p-2}^{(1)}(x)|=1$.

Next we construct the map $Q(X):\B(x,|A_{p-2}(x)|)\to\B(q,|q-1|)$,
by looking at the polygon for $g(x',1+p^{m_0}U)$ expanded in Series B around $u$,
for $x'\in\B(x,|A_{p-2}(x)|)$.
We've already shown $v(d_0(x',u))\equiv v(g(x',q))>0$,
and now claim $v(d_1(x',u))=v(d_{p-2}(x',0))=0$.
Consider the series
$$
d_1(x',u)=\sum_{k=1}kd_k(x',0)u^{k-1}=\sum_{k=1}k\frac{A_k(x')}{(k+2)!}p^{km_0}u^{k-1}
$$
We've seen $v(d_k(x',0))>v(d_{p-2}(x',0))$ unless $k=p^r-2$ for some $r$, and we compute as before,
$v(d_{p^r-2}(x',0))=v(A_{p^r-2}(x))+(p^r-2)\delta-1/(p-1)$, where $\delta=m_0-1/(p-1)$.
Since $A_{p-2}(X)$ divides $A_{p^r-2}(X)$, we have $v(A_{p^r-2}(x))\geq v(A_{p-2}(x))$,
and now it is easy to see $v(d_{p^r-2}(x',0))$ is strictly minimized at $r=1$, hence
$v(d_1(x',u))=v(d_{p-2}(x',0))$.
Since $v(D_{p-2}(x'))\geq 0$,
we have $v(A_{p-2}(x')-A_{p-2}(x))\geq v(x'-x)>v(A_{p-2}(x))$,
hence $v(A_{p-2}(x'))=v(A_{p-2}(x))=1-(p-2)m_0$, and we compute $v(d_{p-2}(x',0))=0$.
Therefore $v(d_1(x',u))=0$.

Since $v(d_0(x',u))>v(d_1(x',u))=0$, the Newton polygon for $g(x',1+p^{m_0}U)$ expanded around $u$ 
shows that for each $x'\in \B(x,|A_{p-2}(x)|)$ there is a root $U=u'\in\B(u,1)$, so that $\bar u'=\bar u$,
and $v(u'-u)=v(d_0(x',u))=v(g(x',q))$.
Setting $q'=1+p^{m_0}u'$, we compute
$$
|q'-q|=p^{-m_0}|g(x',q)|=|([x']_q-x')/x'(x'-1)|
$$
If $|A_{p-2}^{(1)}(x)|=1$, then
$|q'-q|=p^{1-(p-1)m_0}|x'-x|$, so $Q(X)$ is a contraction.
If $|A_{p-2}^{(1)}(x)|<1$, then $|q'-q|<p^{1-(p-1)m_0}|x'-x|$, and $Q(X)$ is not a contraction.
For then the explicit formula for $D_{p-2}(X)$ shows that 
$|D_{p-2}(x')|=1$ in the limit as $|x'-x|$ approaches $1=|A_{p-2}(x)|$,
hence $|q'-q|$ approaches $p^{1-(p-1)m_0}|x'-x|$ arbitrarily closely for $x'\in\B(x,1)$.

This $u'$ is unique since the $p-2$ solutions for $U$ in $g(x',1+p^{m_0}U)$ have
distinct residues, 
by the preceding argument.
Thus we have a well defined map 
$$
Q(X):\B(x,|A_{p-2}(x)|)\to\B(q,|q-1|)
$$ 
sending $x'$ to $q'=Q(x')=1+p^{m_0}u'$.
This map is analytic by Proposition 2.

Next we show $Q(X)$ is surjective.
Suppose $(x,q)\in M$, and $q'\in\B(q,|q-1|)$.
Since $|q'-q|<|q-1|,$ $|q'-1|=|q-1|=p^{-m_0}$.
Define $u$ and $u'$ by $q=1+p^{m_0}u$ and $q'=1+p^{m_0}u'$, and
let $\epsilon=u'-u$, 
so that $v(\epsilon)>0$.
Since $g(x,q)=0$, in Series B we have
$$
g(x,q')=g(x,q')-g(x,q)=\sum_{n=1}^\infty d_n(x,0)((u+\epsilon)^n-u^n)
$$ 
By the binomial theorem, $\epsilon$ is a factor of each $(u+\epsilon)^n-u^n$,
and so $v(g(x,q'))\geq v(\epsilon)>0$.

We show there exists an $x'\in\B(x,|A_{p-2}(x)|)\cap M_{q'}$.
First assume $\bar x\neq 0,1$, so that $v(f(x,q'))=v(g(x,q'))+m_0$.
Let $c_n(x,q')$ be the coefficient of Series A for $f(X,q')$ expanded around $x$.
Since $v(g(x,q'))>0$, we have
$v(c_0(x,q'))=v(g(x,q'))+m_0>m_0\geq (p-1)m_0-1=v(c_p(x,q'))$, 
hence $v(c_0(x,q'))>v(c_p(x,q'))$.
Therefore $f(X,q')$ has a solution $X=x'$ such that $|x'-x|<1$, and 
$x'\in\B(x,1)\cap M_{q'}$.
If $m_0=1/(p-2)$, then $|A_{p-2}(x)|=1$, and we are done.
If $m_0<1/(p-2)$, then the residue multiplicity of $\bar x'$ in $\edge M_{q'}$
is one by Proposition 5, so that $v(c_1(x,q'))=v(c_p(x,q'))=(p-1)m_0-1$,
and by the above computation we have $v(c_0(x,q'))>m_0$.
Since $v(A_{p-2}(x))=1-(p-2)m_0$ by Proposition 4,
we conclude $v(x'-x)=v(c_0(x,q'))-v(c_1(x,q'))>v(A_{p-2}(x))$, as desired.

Next, suppose $\bar x=0$, then $m_0=1/(p-2)$ by Proposition 5,
and $|A_{p-2}(x)|=1$.
Since $0\in\B(x,1)$, and $Q(X):\B(x,1)\to\B(q,|q-1|)$, we may assume $x=0\in M_{q-1}$.
Series A for $f(X,q')$ expanded around $0$
has the trivial fixed point $0$, so $c_0(0,q')=0$/ we will show $v(c_1(0,q'))>v(c_2(0,q'))$.
Since $m_0=1/(p-2)$, we compute $v(c_2(0,q'))=m_0=v(c_p(0,q')$ by $(*)$,
so we have to show $v(c_1(0,q'))>m_0$.
The series for $c_1(0,q')=\log q'/(q'-1)-1$ is
$$
\align
-1+\sum_{n=1}\frac{(-1)^{n+1}}{n}(q'-1)^{n-1}&=\sum_{n=0}\frac{(-1)^{n+1}}{n+2}(q'-1)^{n+1}\\
&=-(q'-1)/2+(q'-1)^2/3-\cdots+(q'-1)^{p-1}/p-\cdots
\endalign
$$
We see at once that only the $n=0$ and $n=p-2$ terms in the right parentheses have
value $m_0=v(q'-1)$.
Since $c_1(0,q)=0$, it suffices to show $v(q'-q)>m_0$ and $v((q'-1)^{p-1}-(q-1)^{p-1})>v(p)+m_0$.
Writing $q'-1=q-1+p^{m_0}\epsilon$ with $v(\epsilon)>0$, as before, we immediately
verify $v(q'-q)>m_0$, and applying
the binomial theorem, we see $v((q-1+p^{m_0}\epsilon)^{p-1}-(q-1)^{p-1})>1+m_0$.
We conclude $v(c_1(0,q'))>v(c_2(0,q'))$, and the Newton polygon shows there exists
an $x'\in\B(0,1)\cap M_{q'}$, as desired.

Next, suppose $\bar x=1$, then again $m_0=1/(p-2)$ by Proposition 5, and $|A_{p-2}(x)|=1$.
Again we may assume $x=1$, since $1\in\B(x,1)$, and the proof 
that $\B(x,1)\cap M_{q'}$ is nonempty is exactly like the $\bar x=0$ case.
We reduce immediately to showing $v(c_1(1,q'))>v(c_2(1,q'))=m_0$.
The series for $c_1(1,q')=q'\log q'/(q'-1)-1$ is
$$
\align
-1+\sum_{n=1}\frac{(-1)^{n+1}}{n}q'(q'-1)^{n-1}&=-1+q'+\sum_{n=1}\frac{(-1)^n}{n+1}q'(q'-1)^n\\
&=-(\frac{q'}{2}-1)(q'-1)+\frac{q'}{3}(q'-1)^2-\cdots+\frac{q'}{p}(q'-1)^{p-1}-\cdots
\endalign
$$
Again
$v(c_1(1,q')>m_0$ is equivalent to $v(c_1(1,q')-c_1(1,q))>m_0$, which is equivalent to
$v(-(q'/2-1)(q'-1)+(q/2-1)(q-1)+q'/p(q'-1)^{p-1}-q/p(q-1)^{p-1})>m_0$,
and the verification is routine.
It follows that $x'\in\B(1,1)\cap M_{q'}$ exists, as desired.

We have shown that for each $(x,q)\in M$ and $q'\in\B(q,|q-1|)$,
there exists an $x'\in\B(x,|A_{p-2}(x)|)\cap M_{q'}$.
Thus the map $Q(X):\B(x,|A_{p-2}(x)|)\to\B(q,|q-1|)$ is onto.
This completes the proof.

\endpf

\flushpar
{\bf Remark.}
If the multiplicity of the residue $\bar x$ in $\edge M_q$ is not equal to one,
then by Proposition 5 and Theorem 6, its multiplicity is two, $\Card(\edge M_q)=p-3$,
$|A_{p-2}(x)|=1$, and $|A_{p-2}^{(1)}(x)|<1$.
It follows that $|A_{p-2}(x')|=1$ and $|A_{p-2}^{(1)}(x')|<1$ for each 
$x'\in\B(x,1)$, and each residue $\bar x'$ has multiplicity two in $\edge M_{Q(x')}$.
We leave aside the problem of proving the existence of such points, and more especially
of proving the existence of nontrivial fixed points $x$ of multiplicity two (in $M_q$).
Of course, if $p=3$, there is no issue.

We next restrict Theorem 6 to the ordinary
$p$-adic integers, which served as the initial motivation for this investigation.

\proclaim{Corollary 7}
Let $M(\Z_p)=\{(x,q)\in M:x,q\in\Z_p\}$.  Then
$$M(\Z_p)\neq\varnothing\;\Longleftrightarrow\;p=3$$
The elements $x\in\Z_3$ that are nontrivial fixed points for some $[X]_q$
form the union $\phi_1(M(\Z_3))=\B(1,1)\cup\B(0,1)$, and
we have an analytic bijection
$$
\align
Q(X):\B(1,1)&\;\longrightarrow\;\B(4,3^{-1})\\
\B(0,1)&\;\longrightarrow\;\B(7,3^{-1})
\endalign
$$
with $|Q(x')-q|=|x'-x|/3$.
\endproclaim

\Pf
By Proposition 3, $M\neq\varnothing$ if and only if $1/(p-1)<v(q-1)\leq 1/(p-2)$,
so we have the first statement.
Assume $p=3$.  
By Proposition 4, $\phi_1(M)\cap\Z_3=\{x\in\Z_3:x\neq 2(\mod 3)\}$,
and by Proposition 3, $\phi_2(M)\cap\Z_3=\{q-1:|q-1|=3^{-1}\}$.
Since the Weierstrass polynomial for Series B has degree one, we see that
$x\in\Z_3$ if and only if $Q(x)=q\in\Z_3$, so these sets are $\phi_1(M(\Z_3))$ and 
$\phi_2(M(\Z_3))$, respectively.
Locally the map $Q(X)$ takes $\B(x,1)$ onto $\B(q,3^{-1})$,
and is a contraction by $1/3$, by Theorem 6.
By sheer luck we find the nontrivial fixed point $x=-1/2$ for $q=4$,
and since $-1/2$ has residue $1$, we conclude that $Q(X)$ takes
$\B(1,1)$ onto $\B(4,3^{-1})$ and $\B(0,1)$ onto $\B(7,3^{-1})$.

\endpf

\flushpar
{\bf Remark.}
Using the $3$-adic Weierstrass preparation theorem together with Series A one can approximate the
nontrivial fixed point of $[X]_q$ for any $q\in\B(4,3^{-1})$ (or $q\in\B(7,3^{-1})$)
to arbitrary accuracy, and conversely using Series B one can approximate the value $q$
for which any $x\in\B(1,1)$ (or $x\in \B(0,1)$) is a fixed point for $[X]_q$.
For example, we find that $x=0$ is a nontrivial fixed point for $[X]_q$, where
$$
q\approx 1+2\cdot 3+3^2+2\cdot 3^3+3^6+2\cdot 3^7
$$

The fixed points of the maps $[X]_q(\mod 3^n):\Z/3^n\to\Z/3^n$ for $q$ such that $|q-1|=3^{-1}$ exhibit a remarkable, seemingly erratic pattern that is nevertheless completely governed by the unique nontrivial fixed point
$x\in\Z_3$.  For example, it can be shown that if $n>2v(x(x-1))+1$, then there are 
$2\cdot 3^{v(x(x-1))+1}+3$ fixed points, and otherwise there are $3^{n-\lfloor n/2\rfloor]}+3$.

As was pointed out earlier,
the nontrivial $3$-adic fixed point for $[X]_4$ is $-1/2=1+3+3^2+3^3+3^4+\cdots$.
In the following table we list
the values of $[X]_4(\mod 3^4)$ on $\Z/3^4\Z$ from $0$ to $3^3$,
which is where the fixed point pattern repeats.  We box the fixed points:
$$
\matrix
\boxed 0 & \boxed 1 & 5 & 21 & \boxed 4 & 17 & 69 & 34 &
56 & 63 & \boxed{10} & 41 & 3 & \boxed{13} & 53 & 51 & 43 & 11\\
45 & \boxed{19} & 77 & 66 & \boxed{22} & 8 & 33 & 52 & 47 & \boxed{27} &
\boxed{28} & 32 & \cdots \\
\endmatrix
$$
By contrast, here is the table for $[X]_4(\mod 3^5)$ on $\Z/3^5\Z$, listed up to $3^4$:
$$
\matrix
\boxed{\ssize 0} & \boxed{\ssize 1} &\ssize 5 &\ssize 21 &\ssize 85 &\ssize 98 &\ssize 150
&\ssize 115 &\ssize
218 &\ssize 144 &\ssize 91 &\ssize 122 &\ssize 3 & \boxed{\ssize 13} &\ssize 53 &\ssize 213 &\ssize 124 &\ssize 11\\
\ssize 45 &\ssize 181 &\ssize 239 &\ssize 228 &\ssize 184 &\ssize 8 &\ssize 33 &\ssize 133 &\ssize 47 &\ssize 189
&\ssize
\boxed{\ssize 28} &\ssize 113 &\ssize 210 &\ssize 112 &\ssize 206 &\ssize 96 &\ssize 142 &\ssize 83 \\
\ssize 90 &\ssize 118 &\ssize 230 &\ssize 192 &\ssize \boxed{\ssize 40} &\ssize 161 &\ssize 159 &\ssize 151 &\ssize 119 &\ssize 234 &\ssize
208 &\ssize 104 &\ssize 174 &\ssize 211 &\ssize 116 &\ssize 222 &\ssize 160 &\ssize 155 \\
\ssize 135 & \boxed{\ssize 55} &\ssize 221 &\ssize 156 &\ssize 189 &\ssize 71 &\ssize 42 &\ssize 169 &\ssize 191 &\ssize 36 &\ssize
145 &\ssize 95 &\ssize 138 & \boxed{\ssize 67} &\ssize 26 &\ssize 105 &\ssize 178 &\ssize 227 \\
\ssize 180 &\ssize 235 &\ssize
212 &\ssize 120 &\ssize 238 &\ssize 224 &\ssize 168 &\ssize 187 &\ssize 20 &\ssize \boxed{\ssize 81} &\ssize \boxed{\ssize 82} &\ssize\cdots\\
\endmatrix
$$

\

\flushpar
{\bf References.}

\roster
\item"{\cite{A}}"
Arens, R.:
{\it Homeomorphisms preserving measure in a group},
Ann. of Math., {\bf 60}, no. 3, (1954), pp. 454--457.
\item"{\cite{B}}"
Bishop, E.:
{\it Isometries of the $p$-adic numbers},
J. Ramanujan Math. Soc. {\bf 8} (1993), no. 1-2, 1--5.
\item"{\cite{C}}"
Conrad, K.:
{\it A $q$-analogue of Mahler expansions. I,}
Adv. Math.  {\bf 153}  (2000),  no. 2, 185--230.
\item"{\cite{D}}"
Dieudonn\'e, J.:
{\it  Sur les fonctions continues $p$-adiques},
Bull. Sci. Math. (2) {\bf 68} (1944), 79--95
\item"{\cite{F}}"
Fray, R. D.:
{\it Congruence properties of ordinary and $q$-binomial coefficients,}
Duke Math. J. {\bf 34} (1967), 467--480.
\item"{[G]}"
F. Q. Gouv\^ea, {\it $p$-adic Numbers, An Introduction} Second edition, Springer-Verlag, New York, 2003.
\item"{\cite{J}}"
Jackson, Rev. F. H.:
{\it $q$-difference equations},
Amer. J. Math. {\bf 32} (1910), 305--314.
\item"{\cite{M}}"
Mahler, K.:
{\it An interpolation series for continuous functions of a $p$-adic variable},
J. Reine Angew. Math. {\bf 199} (1958) 23--34.
\item"{[Se]}"
J.-P. Serre, {\it Lie Algebras and Lie Groups}, LNM 1500, Springer-Verlag, New York, 1992.
\item"{\cite{Su}}"
 Sushchanski\u i, V. I.:
{\it Standard subgroups of the isometry group of the metric space of $p$-adic integers},
V\B isnik Ki\"iv. Un\B iiv. Ser. Mat. Mekh.
{\bf 117} no. 30 (1988), 100--107.
\endroster

\bye